\documentclass[a4paper,12pt,final]{amsart}
\usepackage{times,a4wide,mathrsfs,amssymb,amsmath,amsthm,enumerate,xypic,tikzsymbols,dsfont}

\newcommand{\C}{\mathbb{C}}

\newcommand{\QQ}{\mathbb{Q}}
\newcommand{\NN}{\mathbb{N}}
\newcommand{\PP}{\mathbb{P}}

\newcommand{\OO}{\mathcal O}

\newcommand{\YY}{\mathcal Y}

\newcommand{\JJ}{\mathcal J}

\newcommand{\CC}{\mathcal C}

\newcommand{\MM}{\mathcal M}

\newcommand{\wt}{\widetilde}

\newcommand{\one}{\mathds{1}}

\DeclareMathOperator{\ide}{id}

\DeclareMathOperator{\ima}{Im}

\DeclareMathOperator{\Gr}{Gr}

\DeclareMathOperator{\Dec}{Dec}

\newtheorem{theorem}{Theorem}[section]

\newtheorem{lemma}[theorem]{Lemma}

\newtheorem{corollary}[theorem]{Corollary}
\newtheorem{proposition}[theorem]{Proposition}
\newtheorem{conjecture}[theorem]{Conjecture}
\newtheorem{convention}{Conventions}

\newtheorem{nonumbering}{Theorem}

\newtheorem{nonumberingc}{Corollary}

\theoremstyle{definition}
\newtheorem{remark}[theorem]{Remark}
\newtheorem{definition}[theorem]{Definition}

\newtheorem{nonumberingt}{Acknowledgments}

\begin{document}

\author[Robert Laterveer]
{Robert Laterveer}

\address{Institut de Recherche Math\'ematique Avanc\'ee,
CNRS -- Universit\'e 
de Strasbourg,\
7 Rue Ren\'e Des\-car\-tes, 67084 Strasbourg CEDEX,
FRANCE.}
\email{robert.laterveer@math.unistra.fr}

\title{On the motive of codimension 2 linear sections of $\Gr(3,6)$}

\begin{abstract} We consider Fano sevenfolds $Y$ obtained by intersecting the Grassmannian $\Gr(3,6)$ with a codimension 2 linear subspace (with respect to the Pl\"ucker embedding).
We prove that the motive of $Y$ is Kimura finite-dimensional. We also prove the generalized Hodge conjecture for all powers of $Y$.
 \end{abstract}


\maketitle

\section{Introduction}

\noindent
Given a smooth projective variety $Y$ over $\C$, let $A_i(Y):=CH_i(Y)_{\QQ}$ denote the Chow groups of $Y$ (i.e. the groups of $i$-dimensional algebraic cycles on $Y$ with $\QQ$-coefficients, modulo rational equivalence). Let $A_i^{hom}(Y)\subset A_i(Y)$ denote the subgroup of homologically trivial cycles.

The famous Bloch--Beilinson conjectures \cite{Jan}, \cite{Vo} predict that the Hodge level of the cohomology of $Y$ should have an influence on the size of the Chow groups of $Y$. In particular, there is the following conjecture:

\begin{conjecture}\label{gbc} Let $Y$ be a smooth projective variety of Hodge coniveau $\ge c$ (i.e. the Hodge numbers $h^{p,q}(Y)$ vanish provided $p+q\ge 2c$ and $p<c$). Then
  \[ A_i^{hom}(Y)=0\ \ \ \forall\ i< c\ .\]
  \end{conjecture}
  
(This is known as the ``generalized Bloch conjecture''; for motivation and background cf. \cite[Section 1.2]{Vo} or \cite{Jan}.)  

%

Let $\Gr(3,6)$ be the Grassmannian of 3-dimensional linear subspaces of a fixed 6-dimensional vector space.
In this note, we consider smooth complete intersections
  \[ Y=\Gr(3,6)\cap H_1\cap H_2\ \ \ \subset\ \PP^{19} \]
  of the Grassmannian with 2 Pl\"ucker hyperplanes $H_1, H_2$. This $Y$ is a 7-dimensional Fano variety of Hodge coniveau $3$ (i.e. $h^{p,7-p}(Y)=0$ for $p<3$). 
  The Hodge theory of $Y$ has been studied by Donagi in his thesis \cite{Don}. The derived category of $Y$ has been studied by Deliu \cite{Del} (this derived category is not yet well-understood, because HPD for $\Gr(3,6)$ is still conjectural, cf. Remarks \ref{hpd} and \ref{beyond} below).
  
  Our main result is that Conjecture \ref{gbc} is verified in this case:
  
 \begin{nonumbering}[=Theorem \ref{main}] Let 
  \[ Y:= \Gr(3,6)\cap H_1\cap H_2 \ \ \ \subset\ \PP^{19} \]
  be a smooth dimensionally transverse intersection, where the $H_j$ are Pl\"ucker hyperplanes. Then
    \[ A_i^{hom}(Y)=0\ \ \ \forall\ i\not=3\ .\]
  In particular, $Y$ has finite-dimensional motive (in the sense of \cite{Kim}).
        \end{nonumbering}
        
  The argument proving Theorem \ref{main} uses instances of the {\em Franchetta property} (cf. subsection \ref{ss:fr} below). This is similar to, and inspired by, the seminal work of Voisin on Conjecture \ref{gbc} \cite{V0}, \cite{V1}. 
  
  Theorem \ref{main} has the following consequence:
  
  \begin{nonumberingc}[=Corollary \ref{ghc}] Let $Y$ be as in Theorem \ref{main}. The generalized Hodge conjecture is true for all powers of $Y$.
  \end{nonumberingc}
  
The argument of Corollary \ref{ghc} is as follows: there is a certain elliptic curve $C$ naturally associated to $Y$ (this $C$ is called the {\em Segre curve\/} in honour of C. Segre who studied this curve more than a century ago \cite{Segre}). An equivalent formulation of Theorem \ref{main} is the relation of Chow motives
  \[  h(Y)\ \cong\ h(C)(-3)\oplus \bigoplus\one(\ast)\ \ \ \hbox{in}\ \MM_{\rm rat}\ .\]
  Thus, to prove the generalized Hodge conjecture for powers of $Y$ one is reduced to powers of $C$, for which the generalized Hodge conjecture is known. 
  
  Another application of Theorem \ref{main} concerns Voevodsky's conjecture on smash-equivalence (Corollary \ref{smash}).

It would be interesting to understand more generally the Chow groups of linear sections (of codimension $r>2$) of $\Gr(3,6)$, cf. Remark \ref{beyond}.

 \vskip0.8cm

\begin{convention} In this note, the word {\sl variety\/} will refer to a reduced irreducible scheme of finite type over $\C$. A {\sl subvariety\/} is a (possibly reducible) reduced subscheme which is equidimensional. 

{\bf All Chow groups will be with rational coefficients}: we denote by $A_j(Y):=CH_j(Y)_{\QQ} $ the Chow group of $j$-dimensional cycles on $Y$ with $\QQ$-coefficients; for $Y$ smooth of dimension $n$ the notations $A_j(Y)$ and $A^{n-j}(Y)$ are used interchangeably. 
The notations $A^j_{hom}(Y)$ and $A^j_{AJ}(X)$ will be used to indicate the subgroup of homologically trivial (resp. Abel--Jacobi trivial) cycles.

The contravariant category of Chow motives (i.e., pure motives with respect to rational equivalence as in \cite{Sc}, \cite{MNP}) will be denoted 
$\MM_{\rm rat}$.
\end{convention}

 \vskip1.0cm

\section{Preliminaries}

\subsection{Codimension 1 linear sections} 

As a warm-up for the codimension 2 case, let us first consider codimension 1 linear sections of $\Gr(3,6)$.

\begin{theorem}[Donagi \cite{Don}]\label{don}  Let 
  \[ Y:= \Gr(3,6)\cap H \ \ \ \subset\ \PP^{19} \]
  be a smooth dimensionally transverse intersection, where $H$ is a Pl\"ucker hyperplane.
  The interesting Hodge numbers of $Y$ are
    \[ h^{p,8-p}(Y)=\begin{cases}  4 & \hbox{if}\ p=4\ ,\\
                                                      0 & \hbox{otherwise}\ .
                                                      \end{cases}\]
                \end{theorem}
             
        \begin{proof} This is contained in \cite[Section 3.4]{Don}.
          \end{proof}      
          
    \begin{remark}\label{element} Let $Y= \Gr(3,6)\cap H$ be a hypersurface as in Theorem \ref{don}. As the Hodge coniveau of $Y$ is $4$, \, Conjecture \ref{gbc} 
    (combined with the Bloch--Srinivas argument \cite{BS}, \cite{moi}) predicts that 
      \begin{equation}\label{predict} A_i^{hom}(Y)=0\ \ \ \forall \ i\ .\end{equation}
   In this case, this is readily verified using a construction of Donagi's \cite{Don}: we choose a hyperplane $P_5\subset V_6$, and we consider the incidence variety
     \[ \wt{\Gr}  := \Bigl\{ (A,\ell)\in \Gr(3,V_6)\times \Gr(2,P_5)\ \Big\vert\, \ell\subset A\Bigr\}\ \ \ \subset\ \Gr(3,V_6)\times \Gr(2,P_5)\ .\]
    The first projection $\wt{\Gr}\to \Gr(3,V_6)$ is birational (it is actually a blow-up with center the locus $\sigma_{111}(P)$ of subspaces contained in $P$). The second projection
    \[ \Pi\colon\ \ \wt{\Gr}\ \to\ \Gr(2,P_5)\]
    is a $\PP^3$-fibration.
   We now consider the morphism
     \[ \Pi_Y\colon\ \ \wt{Y}\ \to\ \Gr(2,P_5) \ ,\]
     obtained by restricting $\Pi$ to the strict transform $\wt{Y}$ of $Y$ in $\wt{\Gr}$. As explained in \cite[Section 3.4]{Don}, for $P_5$ generic with respect to $Y$, the morphism
     $\Pi_Y$ is a $\PP^2$-fibration over $\Gr(2,P_5)\setminus S$, and a $\PP^3$-fibration over $S$, where $S\subset\Gr(2,P_5)$ is a closed subvariety isomorphic to $\PP^1\times\PP^1$. Applying the Chow-theoretic version of Cayley's trick \cite[Theorem 3.1]{Jiang} to this set-up, we find that $A_i(\wt{Y})$ is a direct sum of Chow groups of $\Gr(2,5)$ and of $S$, hence $\wt{Y}$ has trivial Chow groups. It follows that $Y$ also has trivial Chow groups, i.e. the prediction \eqref{predict} is verified. (Another proof of \eqref{predict} is given in \cite[Theorem 3.2]{hypgr}.)
     \end{remark}                                                    

\begin{remark}\label{hpd0} Let $Y= \Gr(3,6)\cap H$ be a hypersurface as in Theorem \ref{don}. The theory of homological projective duality \cite{Kuz0}, \cite{Kuz1}, \cite{Thom} suggests that the derived category of $Y$ should admit a full exceptional collection. The construction of the HPD dual of $\Gr(3,6)$ appears to be an open problem (cf. \cite[Conjecture 21]{Del}, where a non-commutative resolution of the double cover of $\PP^{19}$ branched along a certain quartic hypersurface is suggested as a candidate). Nevertheless, it seems likely the existence of a full exceptional collection for $Y$ can be proven by looking at the above construction and applying a categorical version of Cayley's trick (\cite[Theorem 2.6]{JL} or \cite[Proposition 47]{BFM}).
\end{remark}

\subsection{Codimension 2 linear sections}

\begin{theorem}[Donagi \cite{Don}]\label{don2}  Let 
  \[ Y:= \Gr(3,6)\cap H_1\cap H_2 \ \ \ \subset\ \PP^{19} \]
  be a smooth dimensionally transverse intersection, where the $H_j$ are Pl\"ucker hyperplanes.
  The Hodge diamond of $Y$ is
     \[ \begin{array}[c]{ccccccccccccccc}
     && &&&&& 1 &&&&&&&\\
     && &&&&0&&0&&&&&&\\
     && &&&0&&1&&0&&&&&\\
    &&  &&0&&0&&0&&0&&&&\\
     && &0&&0&&2&&0&&0&&&\\
       && 0 &&0&&0&&0&&0&&0&&\\     
      &0& &0&&0&&3&&0&&0&&0&\\
      0&&  0&&0&&1&&1&&0&&0&&0\\
           &0& &0&&0&&3&&0&&0&&0&\\
              && 0 &&0&&0&&0&&0&&0&&\\              
            && &0&&0&&2&&0&&0&&&\\
           &&  &&0&&0&&0&&0&&&&\\
    && &&&0&&1&&0&&&&&\\
       && &&&&0&&0&&&&&&\\
 && &&&&& 1 &&&&&&&\\
            \end{array}\]

             \end{theorem}
             
        \begin{proof} This is contained in \cite[Section 3.4]{Don}. (Alternatively, one could apply \cite{FM0} to find the Hodge number $h^{4,3}(Y)$.)
          \end{proof}      

\begin{corollary}\label{ghc1} Let $Y:=\Gr(3,6)\cap H_1\cap H_2$ be as in Theorem \ref{don2}. There exist an elliptic curve $C$ and a correspondence $\Gamma\in A^{4}(C\times Y)$ inducing an isomorphism
  \[  \Gamma_\ast\colon\ \ H^1(C,\QQ)\ \xrightarrow{\cong}\ H^7(Y,\QQ)\ .\]
  \end{corollary}
  
 \begin{proof} Let $\JJ^4(Y)$ denote the intermediate Jacobian. Because the Hodge coniveau of $H^7(Y,\QQ)$ is 3, $\JJ^4(Y)$ is an abelian variety. Because $h^{4,3}(Y)=1$, the dimension of $\JJ^4(Y)$ is 1, i.e. $\JJ^4(Y)$ is an elliptic curve. As explained for instance in \cite[Proof of Lemma 3]{Murre}, the general theory of abelian varieties guarantees the existence of an elliptic curve $C$ and a correspondence
  $\Gamma\in A^{4}(C\times Y)$ such that there is a commutative diagram
    \[ \begin{array}[c]{ccc}
       A^1_{hom}(C) &  \xrightarrow{\Gamma_\ast} & A^4_{hom}(Y)\\
       &&\\
      \ \  \downarrow{\scriptstyle AJ}&&   \ \  \downarrow{\scriptstyle AJ}\\
      &&\\
     C\cong \JJ^1(C) & \xrightarrow{\Gamma_\ast}& \ \, \JJ^4(Y)\ .\\
      \end{array}\]
      Here $AJ$ is the Abel--Jacobi map, and $\JJ^1(C)$ is the Jacobian of $C$. The left arrow is an isomorphism, and the lower horizontal arrow is an isogeny of elliptic curves.
      It follows that $\Gamma$ induces an isomorphism
      \[ \Gamma_\ast\colon\ \ H^{0,1}(C)\ \xrightarrow{\cong}\ H^{3,4}(Y)\ ,\]
      and hence (taking the complex conjugate) also
      \[ \Gamma_\ast\colon\ \ H^{1,0}(C)\ \xrightarrow{\cong}\ H^{4,3}(Y)\ .\]
     Invoking the Hodge decomposition $H^7(Y,\C)=H^{4,3}(Y)\oplus H^{3,4}(Y)$ (and likewise $H^1(C,\C)=H^{1,0}(C)\oplus H^{0,1}(C)$), we find that $\Gamma$ induces an isomorphism
      \[  \Gamma_\ast\colon\ \ H^1(C,\C)\ \xrightarrow{\cong}\ H^7(Y,\C)\ .\]            
\end{proof}

\begin{remark} It is shown by Donagi \cite{Don} that one can actually find a geometric incarnation of the elliptic curve $C$ of Corollary \ref{ghc1}. There is a certain elliptic curve in $\Gr(3,6)$ naturally associated to the pencil of hyperplanes defining $Y$, this is called the {\em Segre curve\/} in \cite[Section 3.3]{Don} (in honour of C. Segre who had already studied this curve \cite{Segre}). Donagi proves \cite[Theorem 3.8]{Don} that there is a natural isomorphism (with a geometric interpretation) from the Segre elliptic curve to the intermediate Jacobian $\JJ^4(Y)$.
\end{remark}

\begin{remark} Let $Y:=\Gr(3,6)\cap H_1\cap H_2$ be as in Theorem \ref{don2}. To understand the Chow groups of $Y$, it is natural to try and apply the method of Remark \ref{element}.
That is, one would like to consider the morphism
  \[ \Pi_Y\colon\ \ \wt{Y}\ \to\ \Gr(2,5) \]
  obtained by restricting the morphism $\Pi\colon\wt{\Gr}\to \Gr(2,5)$ (of Remark \ref{element}) to the strict transform $\wt{Y}$ of $Y$ in $\wt{\Gr}$. Unfortunately, this approach runs into trouble: the locus
  $S\subset\Gr(2,5)$ (where the fibers of $\Pi_Y$ have larger dimension) is now a quadric surface bundle with some singular fibers, and it seems difficult to handle the Chow groups of $S$.

For this reason, I have preferred to use the ``spread'' method (which means considering the universal family of sections $Y$) to prove Theorem \ref{main}.
\end{remark}

\begin{remark}\label{hpd} Let $Y:=\Gr(3,6)\cap H_1\cap H_2$ be as in Theorem \ref{don2}. Homological projective duality predicts that there is a semi-orthogonal decomposition
  \[ D^b(Y)=\bigl\langle D^b(C), A_1, \ldots, A_r\bigr\rangle\ ,\]
  where $C$ is the Segre elliptic curve and the $A_j$ are exceptional objects. As mentioned in Remark \ref{hpd0}, a conjectural candidate for the HPD dual of $\Gr(3,6)$ is identified in \cite[Conjecture 21]{Del}.
  \end{remark}

  \subsection{The Franchetta property}
 \label{ss:fr}

 \begin{definition} Let $\YY\to B$ be a smooth projective morphism, where $\YY, B$ are smooth quasi-projective varieties. We say that $\YY\to B$ has the {\em Franchetta property in codimension $j$\/} if the following holds: for every $\Gamma\in A^j(\YY)$ such that the restriction $\Gamma\vert_{Y_b}$ is homologically trivial for the very general $b\in B$, the restriction $\Gamma\vert_{Y_b}$ is rationally trivial, i.e. $\Gamma\vert_{Y_b}$ is zero in $A^j(Y_b)$ for all $b\in B$.
 
 We say that $\YY\to B$ has the {\em Franchetta property\/} if $\YY\to B$ has the Franchetta property in codimension $j$ for all $j$.
 \end{definition}
 
 This property is studied in \cite{PSY}, \cite{BL}, \cite{FLV}, \cite{FLV3}.
 
 \begin{definition} Given a family $\YY\to B$ as above, with $Y:=Y_b$ a fiber, we write
   \[ GDA^j_B(Y):=\ima\Bigl( A^j(\YY)\to A^j(Y)\Bigr) \]
   for the subgroup of {\em generically defined cycles}. 
  In a context where it is clear to which family we are referring, the index $B$ will often be suppressed from the notation.
  \end{definition}
  
  With this notation, the Franchetta property amounts to saying that $GDA^\ast_B(Y)$ injects into cohomology, under the cycle class map. 
 
   \subsection{A Franchetta-type result}

  \begin{proposition}\label{Frtype} Let $M$ be a smooth projective variety with 
  $A^\ast_{hom}(M)=0$. Let $L_1,\ldots,L_r\to {M}$ be very ample line bundles, and let
  $\YY\to B$ be the universal family of smooth dimensionally transverse complete intersections of type 
    \[ Y={M}\cap H_1\cap\cdots\cap H_r\ ,\ \ \  H_j\in\vert L_j\vert\ .\]
  Assume the fibers $Y=Y_b$ have $H^{\dim Y}_{tr}(Y,\QQ)\not=0$.
  There is an inclusion
    \[ \ker \Bigl( GDA^{\dim Y}_B(Y\times Y)\to H^{2\dim Y}(Y\times Y,\QQ)\Bigr)\ \ \subset\ \Bigl\langle (p_1)^\ast GDA^\ast_B(Y), (p_2)^\ast GDA^\ast_B(Y)  \Bigr\rangle\ ,\]
  where $p_1, p_2$ denote the projection from $Y\times Y$ to first resp. second factor.
   \end{proposition}
   
   \begin{proof} This is essentially Voisin's ``spread'' result \cite[Proposition 1.6]{V1} (cf. also \cite[Proposition 5.1]{LNP} for a reformulation of Voisin's result). We give a different proof based on \cite{FLV}. Let $\bar{B}:=\PP H^0({M},L_1\oplus\cdots\oplus L_r)$ (so $B\subset \bar{B}$ is a Zariski open), and let us consider the projection
   \[ \pi\colon\ \  \YY\times_{\bar{B}} \YY\ \to\ M\times M\ .\]
   Using the very ampleness assumption, one finds that $\pi$ is a $\PP^s$-bundle over $({M}\times {M})\setminus \Delta_{{M}}$, and a $\PP^t$-bundle over the diagonal $\Delta_{{M}}$.
   That is, $\pi$ is what is termed a {\em stratified projective bundle\/} in \cite{FLV}. As such, \cite[Proposition 5.2]{FLV} implies the equality
      \begin{equation}\label{stra} GDA^\ast_B(Y\times Y)= \ima\Bigl( A^\ast(M\times M)\to A^\ast(Y\times Y)\Bigr) +  \Delta_\ast GDA^\ast_B(Y)\ ,\end{equation}
      where $\Delta\colon Y\to Y\times Y$ is the inclusion along the diagonal. The assumption $A^\ast_{hom}(M)=0$ implies that $M$ has the Chow--K\"unneth property, i.e. $A^\ast(M\times M)$ is isomorphic to $A^\ast(M)\otimes A^\ast(M)$ (this follows from \cite[Proposition 4.22]{Vo}). 
      Base-point freeness of the $L_j$ implies that 
        \[  GDA^\ast_B(Y)=\ima\bigl( A^\ast(M)\to A^\ast(Y)\bigr)\ .\]
       The equality \eqref{stra} thus reduces to
      \[ GDA^\ast_B(Y\times Y)=\Bigl\langle (p_1)^\ast GDA^\ast_B(Y), (p_2)^\ast GDA^\ast_B(Y), \Delta_Y\Bigr\rangle\ \]   
      (where $p_1, p_2$ denote the projection from $S\times S$ to first resp. second factor). The assumption that $Y$ has non-zero transcendental cohomology
      implies that the class of the diagonal $\Delta_Y$ is not decomposable in cohomology (indeed, decomposable correspondences act as zero on the transcendental cohomology). It follows that
      \[ \begin{split}  \ima \Bigl( GDA^{\dim Y}_B(Y\times Y)\to H^{2\dim Y}(Y\times Y,\QQ)\Bigr) =&\\
       \ima\Bigl(  \Dec^{\dim Y}(Y\times Y)\to H^{2\dim Y}(Y\times Y,\QQ)\Bigr)& \oplus \QQ[\Delta_Y]\ ,\\
       \end{split}\]
      where we use the shorthand 
       \[ \Dec^j(Y\times Y):= \Bigl\langle (p_1)^\ast GDA^\ast_B(Y), (p_2)^\ast GDA^\ast_B(Y)\Bigr\rangle\cap A^j(Y\times Y) \ \]     
       for the {\em decomposable cycles\/}. 
       We now see that if $\Gamma\in GDA^{\dim Y}(Y\times Y)$ is homologically trivial, then $\Gamma$ does not involve the diagonal and so $\Gamma\in \Dec^{\dim Y}(Y\times Y)$.
       This proves the proposition.
         \end{proof}
  
  \begin{remark} Proposition \ref{Frtype} has the following consequence: if the family $\YY\to B$ has the Franchetta property, then $\YY\times_B \YY\to B$ has the Franchetta property in codimension $\dim Y$.
   \end{remark}
 
 \subsection{A Chow--K\"unneth decomposition}

\begin{lemma}\label{ck} Let $M$ be a smooth projective variety with $A^\ast_{hom}(Y)=0$. Let $Y\subset M$ be a smooth complete intersection as in Proposition \ref{Frtype}, of dimension $\dim Y=d$.
 The variety $Y$ has a self-dual Chow--K\"unneth decomposition $\{\pi^j_Y\}$ with the property that
   \[  h^j(Y):=(Y,\pi^j_Y,0) =\oplus \one(\ast)\ \ \ \hbox{in}\ \MM_{\rm rat}\ \ \ \forall\ j\not=d \ .\]
   
   Moreover, this decomposition is {\em generically defined\/}: writing $\YY\to B$ for the universal family (of complete intersections of the type of $Y$),
 there exist relative projectors $\pi^j_\YY\in A^{d}(\YY\times_B \YY)$ such that $\pi^j_Y=\pi^j_\YY\vert_{b}$ (where $Y=Y_b$ for $b\in B$). 
     \end{lemma} 
 
 \begin{proof} This is a standard construction, one can look for instance at \cite{Pet} (in case $d$ is odd, which will be the case in this note,  the ``variable motive'' $h(Y)^{var}$ of 
 \cite[Theorem 4.4]{Pet}  coincides with $h^{d}(Y)$).
\end{proof}

   \section{Main result}
   
   \begin{theorem}\label{main} Let 
  \[ Y:= \Gr(3,6)\cap H_1\cap H_2 \ \ \ \subset\ \PP^{19   } \]
  be a smooth dimensionally transverse intersection with 2 hyperplanes $H_1, H_2$ (with respect to the Pl\"ucker embedding). Let $C$ be the Segre curve associated to $Y$. There is an isomorphism of motives
    \begin{equation}\label{isom}  h(Y)\ \cong\ h(C)(-3)\oplus \bigoplus\one(\ast)\ \ \ \hbox{in}\ \MM_{\rm rat}\ .\end{equation}
  In particular, $Y$ has finite-dimensional motive (in the sense of \cite{Kim}), and
          \[ A^i_{AJ}(Y)=0\ \ \ \forall i\ .\]
     \end{theorem}
   
   \begin{proof} Let $\YY\to B$ denote the universal family of smooth dimensionally transverse intersections as in the theorem, where $B$ is a Zariski open in
     \[ \bar{B}:= \PP H^0(\Gr(3,6),\OO_{\Gr(3,6)}(1)^{\oplus 2})\ .\]
   Proposition \ref{Frtype} applies to this set-up (with $M=\Gr(3,6)$), and gives an inclusion
       \begin{equation}\label{inc} \ker \Bigl( GDA^{7}_B(Y\times Y)\to H^{14}(Y\times Y,\QQ)\Bigr)\ \ \subset\ \Bigl\langle (p_1)^\ast GDA^\ast_B(Y), (p_2)^\ast GDA^\ast_B(Y) 
        \Bigr\rangle\ .\end{equation}
        
Let us construct an interesting cycle in $GDA^{7}_B(Y\times Y)$ to which we can apply \eqref{inc}. For any $Y=Y_b$ with $b\in B$, Corollary \ref{ghc1} gives us a smooth curve $C=C_b$ and a cycle $\Gamma\in A^{4}(C\times Y)$ inducing a surjection
  \[ \Gamma_\ast\colon\ \     H^1(C,\QQ)\ \twoheadrightarrow\ H^7(Y,\QQ)\ .\]
  Writing $\CC\to B$ for the universal family of Segre curves, the cycle $\Gamma$ naturally exists relatively, i.e. $\Gamma\in GDA^4_B(F\times Y)$. 
Since both $C$ and $Y$ verify the standard conjectures, the right-inverse to $\Gamma_\ast$ is correspondence-induced, i.e. there exists $\Psi\in A^{4}(Y\times C)$ such that
  \[ (\Gamma\circ \Psi)_\ast =\ide\colon\ \ H^7(Y,\QQ)\ \to\ H^7(Y,\QQ) \ \]
  (This follows as in \cite[Proof of Proposition 1.1]{V4}).\footnote{It seems likely that one can actually take $\Psi$ to be a multiple of the transpose of $\Gamma$. For this, one would need to know that the map $\Gamma_\ast$ of Corollary \ref{ghc1} is compatible with cup-product up to a multiple; this is the point of view taken in \cite{V1} to create the inverse correspondence $\Psi$.}
  
  We now involve the (generically defined) Chow--K\"unneth decomposition $\pi^j_Y\in A^7(Y\times Y)$ given by Lemma \ref{ck}. The above means that for $Y=Y_b$ for any $b\in B$, there is vanishing
    \[  (\Delta_Y -\Gamma\circ \Psi)\circ \pi^7_Y=0\ \ \ \hbox{in}\ H^{14}(Y\times Y,\QQ)\ .\]
  
  Applying Voisin's Hilbert schemes argument \cite[Proposition 3.7]{V0}, \cite[Proposition 4.25]{Vo} (cf. also \cite[Proposition 2.10]{Lfam} for the precise form used here), we can assume that $\Psi$ is also generically defined, and hence
    \[  (\Delta_Y -\Gamma\circ \Psi)\circ \pi^7_Y\ \ \in\ GDA^7(Y\times Y)\ .\]
  Now looking at \eqref{inc}, we learn that this cycle is {\em decomposable\/}, i.e.
   \[      (\Delta_Y -\Gamma\circ \Psi)\circ \pi^7_Y\ \ \in\ \Bigl\langle (p_1)^\ast GDA^\ast(Y), (p_2)^\ast GDA^\ast(Y)  \Bigr\rangle\ . \]
  That is, for any $Y=Y_b$ with $b\in B$ we can write
   \[   (\Delta_Y -\Gamma\circ \Psi)\circ \pi^7_Y=\gamma\ \ \ \hbox{in}\ A^7(Y\times Y)\ ,\]   
   with $\gamma\in A^\ast(Y)\otimes A^\ast(Y)$. Since the $\pi^j_Y, j\not=7$ of Lemma \ref{ck} are also decomposable (i.e. they are in $A^\ast(Y)\otimes A^\ast(Y)$), 
   this implies that we can write
   \[   \Delta_Y -\Gamma\circ \Psi=\gamma^\prime\ \ \ \hbox{in}\ A^7(Y\times Y)\ ,\]      
   with $\gamma^\prime\in A^\ast(Y)\otimes A^\ast(Y)$. Being decomposable, $\gamma^\prime$ does not act on Abel--Jacobi trivial cycles, and so
   \[ A^i_{AJ}(Y)\ \xrightarrow{\Psi_\ast}\ A^{i-3}_{AJ}(C)\ \xrightarrow{\Gamma_\ast}\ A^i_{AJ}(Y) \]
   is the identity. But $C$ being a curve, the group in the middle vanishes for all $i$. This proves the vanishing
     \[   A^\ast_{AJ}(Y)=0 \] 
 for all $Y=Y_b$. The Kimura finite-dimensionality of $Y$ also follows, since submotives of sums of motives of curves are finite-dimensional.
 
 We have now proven that there is a split injection
   \[ h^7(Y)\ \hookrightarrow\ h^1(C)(-3)\ \ \ \hbox{in}\ \MM_{\rm rat} \ .\]
   Since the motive $h^1(C)$ is indecomposable and $h^7(Y)$ is non-zero, this injection is actually an isomorphism
   \begin{equation}\label{act} h^7(Y)\ \xrightarrow{\cong}\ h^1(C)(-3)\ \ \ \hbox{in}\ \MM_{\rm rat} \ .\end{equation}  
   Combining \eqref{act} with the equalities
     \[ \begin{split} h(Y) &= h^7(Y)\oplus \one \oplus \one(-1) \oplus \one(-2)^{\oplus 2} \oplus \one(-3)^{\oplus 3} \oplus  \one(-4)^{\oplus 3}\oplus  \one(-5)^{\oplus 2} \oplus \one(-6) \oplus \one(-7),\\  
     h(C)&= h^1(C)\oplus \one \oplus \one(-1)\ \ \ \ \ 
      \ \ \ \hbox{in}\ \MM_{\rm rat}\ \\
      \end{split}\]
      (cf. Lemma \ref{ck}),
    this gives the isomorphism \eqref{isom}.
    \end{proof}

   \section{Two consequences}
   
   \begin{corollary}\label{ghc} Let $Y$ be as in Theorem \ref{main}. Then the generalized Hodge conjecture is true for $Y^m$ for all $m\in\NN$.
   \end{corollary} 
   
   \begin{proof} The isomorphism of motives of Theorem \ref{main} implies there is an isomorphism of Hodge structures
    \[ H^j(Y^m,\QQ)\ \cong\  H^{j-6m}(C^m,\QQ)(-3m) \oplus \bigoplus H^\ast(C^{m-1},\QQ)(\ast)\oplus  \cdots \oplus \bigoplus \QQ(\ast)\ .\]
   Since this isomorphism is also compatible with the coniveau filtration \cite[Proposition 1.2]{V4}, one is reduced to proving the generalized Hodge conjecture for powers of $C$. This is known thanks to work of Abdulali \cite[Section 8.1]{Ab} (cf. also \cite[Corollary 3.13]{Vial}). 
    \end{proof}
   
  \begin{remark} Corollary \ref{ghc} does not really use the full force of Theorem \ref{main}: to prove Corollary \ref{ghc}, it suffices to have an isomorphism of homological motives linking $Y$ and $C$; such an isomorphism follows readily from Corollary \ref{ghc1}.
  \end{remark}
   
For the next consequence, we recall that a cycle $a\in A^i(Y)$ is called {\em smash-nilpotent\/} if $a^{\times N}$ is zero in $A^{Ni}(Y^N)$ for some $N\in\NN$. Two cycles are {\em smash-equivalent\/} if their difference is smash-nilpotent.
 Voevodsky has conjectured that smash-equivalence coincides with numerical equivalence for all smooth projective varieties \cite{Voe}. Using Theorem \ref{main}, we verify this in some cases:  
   
 \begin{corollary}\label{smash}  Let $Y$ be as in Theorem \ref{main}. Then smash-equivalence and numerical equivalence coincide on $Y^m$ for $m\le 3$.
 \end{corollary}
 
 \begin{proof} The isomorphism of motives of Theorem \ref{main} implies there is an isomorphism of Chow groups
   \[ A^j(Y^m)\ \cong\ A^{j-3m}(C^m) \oplus \bigoplus A^\ast(C^{m-1})\oplus  \cdots \oplus \bigoplus \QQ\ ,\]
   respecting any adequate equivalence relation. The result follows, since smash-equivalence and numerical equivalence coincide (for zero-cycles and divisors \cite{Voe} and hence) for all surfaces, and for abelian threefolds \cite{KS}.
    \end{proof}

 \section{And beyond ?}
 
 \begin{remark}\label{beyond} What can one say about the Chow groups (or Chow motive) of smooth dimensionally transverse intersections
   \[ Y:= \Gr(3,6)\cap H_1\cap \cdots \cap H_r\ \ \ \subset\ \PP^{19}\ \]
 for arbitrary $r$ ?
The conjectural HPD picture drawn in \cite[Section 5.3]{Del} suggests the following prediction: for $r\le 5$, there is an isomorphism of motives
  \[  h(Y)\oplus\bigoplus\one(\ast)\ \cong\ h(X)(r-5)\oplus \bigoplus\one(\ast)\ \ \ \hbox{in}\ \MM_{\rm rat}   \ ,\]
  where $X$ is a double cover of $\PP^{r-1}$ branched along a certain smooth quartic hypersurface. (For $r>5$, $X$ should be a certain resolution of singularities of a double cover
  of $\PP^{r-1}$ branched along a singular quartic.)
  
  To prove this isomorphism for $r>2$ using the method of this note, one is essentially reduced to proving this isomorphism holds modulo homological equivalence and is generically defined.
  \end{remark}

 \vskip1cm
\begin{nonumberingt} Thanks to the referee for helpful comments. Thanks to Mama-san of the izakaya in Schiltigheim.
\end{nonumberingt}

\end{document}